\documentclass[10pt]{amsart}

\usepackage[latin1]{inputenc}
\usepackage{amsmath}
\usepackage{amsfonts}
\usepackage{amssymb}
\usepackage{ytableau}
\usepackage[mathscr]{eucal}
\usepackage{diagrams}
\usepackage{xy,color,graphicx}
\usepackage{hyperref}
\xyoption{all}

\newcommand{\wis}[1]{{\text{\em \usefont{OT1}{cmtt}{m}{n} #1}}}

\newcommand{\N}{\mathbb{N}}

\newcommand{\mon}{\mathbb{N}^{\times}_+}

\newtheorem{definition}{Definition}
\newtheorem{proposition}{Proposition}

\newtheorem{theorem}{Theorem}

\title{Covers of the arithmetic site}
\author{Lieven Le Bruyn} 
\address{Department of Mathematics, University of Antwerp \\ 
 Middelheimlaan 1, B-2020 Antwerp (Belgium) \\ {\tt lieven.lebruyn@uantwerpen.be}}

\begin{document}
\sloppy

\maketitle

\begin{abstract} We give an explicit description of the Barr- and Diaconescu covers of the arithmetic site, which are relevant to cohomology. Further, we construct the arithmetic site as the commutative shadow of a noncommutative topological space.
\end{abstract}

\section{The arithmetic site}

The {\em arithmetic site}, introduced and studied by A. Connes and C. Consani in \cite{CC} and \cite{CC2}, is an algebraic geometric object deeply related to the non-commutative geometric approach to the Riemann hypothesis. It involves two elaborate mathematical concepts: the notion of {\em topos} in geometry and that of {\em characteristic $1$ structures} in algebra.

It is defined as a topos $\widehat{\wis{C}}$ endowed with a structure sheaf $\mathcal{O}$. Here, $\widehat{\wis{C}}$ is the {\em topos} of sheaves on the small category $\wis{C}$ consisting of one object $\bullet$ with monoid of endomorphisms isomorphic to the multiplicative semigroup $\mathbb{N}^{\times}_+$ of strictly positive integers and equipped with the {\em chaotic Grothendieck topology}. Thus, $\widehat{\wis{C}}$ is the category of contravariant functors from $\wis{C}$ to $\wis{Sets}$, the category of sets. Any such functor $\wis{C}^{op} \rTo \wis{Sets}$ is fully determined by a set $X = X(\bullet)$, equipped with a commuting family of endomorphisms $\Psi_n : X \rTo X$ for all $n \in \mathbb{N}^{\times}_+$, with $\Psi_1$ the identity morphism.

That is, $\widehat{\wis{C}}$ is equivalent to the category $\wis{Sets}-\mathbb{N}^{\times}_+$ of all sets equipped with a right action by the monoid $\mathbb{N}^{\times}_+$ and action preserving maps as morphisms. The connection with {\em characteristic one structures} comes from the fact that this action is fully determined by the commuting endomorphisms $\Psi_p$ for $p \in \mathbb{P}$ (the set of all prime numbers), which can be viewed as {\em Frobenius morphisms}.

However, there is a large group of automorphisms of the topos $\widehat{\wis{C}}$ arising from automorphisms of the monoid $\wis{C}$ (such as arbitrary permutations of the prime numbers). For this reason, Connes and Consani equip $\widehat{\wis{C}}$ with a {\em structure sheaf}, turning the presheaf topos $\widehat{\wis{C}}$ in the arithmetic site $(\widehat{\wis{C}},\mathbb{Z}_{max})$. Here, $\mathbb{Z}_{max}$ is the fundamental semiring in characteristic $1$. That is, $\mathbb{Z}_{max} = (\mathbb{Z} \cup \{ - \infty \}, max, +)$, which means that ordinary addition is replaced by $x \oplus y = max(x,y)$ and ordinary multiplication by $x \otimes y = x+y$, see \cite{CC} and \cite{CC2}.

Two topos-theoretic objects associated to $\widehat{\wis{C}}$ carry a lot of arithmetic information. The {\em points} of the topos can be identified with the finite ad\`ele classes \cite{CC2} and the {\em subobject-classifier} is a complete Heyting algebra on the set of all (right) ideals of $\mathbb{N}^{\times}_+$.

\subsection{The points of $\widehat{\wis{C}}$ : } Recall that a {\em point $p$ of the topos $\widehat{\wis{C}}$} is a geometric morphism from $\wis{Sets}$ to $\widehat{\wis{C}}$, see \cite[VII.5]{MM}. Using Grothendieck's description via filtering functors, see \cite[VII.6]{MM}, a point $p$ corresponds to a set $P$ having a free left- $\mathbb{N}^{\times}_+$-action such that $P$ is of rank one meaning that for all $x,y \in P$ there exists an element $z \in P$ and numbers $n,m \in \mathbb{N}^{\times}_+$ such that $x=n.z$ and $y=m.z$. As such, points correspond to additive subgroups of $\mathbb{Q}$. If $\mathbb{A}^f$ is the ring of finite ad\`eles of $\mathbb{Q}$, that is $\mathbb{A}^f = \mathbb{Q} \otimes \prod_{p \in \mathbb{P}} \widehat{\mathbb{Z}}_p$, then it is proved in \cite[Prop. 2.5]{CC2} that any such additive subgroup is of the form
\[
P_a = \{ q \in \mathbb{Q}~|~a.q \in \widehat{\mathbb{Z}} \} \quad \text{for some $a \in \mathbb{A}^f/\widehat{\mathbb{Z}}^*$} \]
where $\widehat{\mathbb{Z}}^*$ is the multiplicative group of invertible elements in the ring of profinite integers $\widehat{\mathbb{Z}} = \prod_{p \in \mathbb{P}} \widehat{\mathbb{Z}}_p$ acting by (right) multiplication on $\mathbb{A}^f$. Moreover, the map $a \mapsto P_a$ induces a canonical bijection
\[
\wis{pts}(\widehat{\wis{C}}) = \mathbb{Q}^{\times}_+ \backslash \mathbb{A}^f / \widehat{\mathbb{Z}}^* = [ \mathbb{A}^f ]  \]
between the isomorphism classes of points of the arithmetic site and finite ad\`ele classes $[\mathbb{A}^f]$. see \cite[Prop. 2.5]{CC2}. Alternatively one can describe additive subgroups of $\mathbb{Q}$ by supernatural numbers as in \cite{BZ}. Recall that a {\em supernatural number} (also called a is a formal product $s = \prod_{p \in \mathbb{P}} p^{s_p}$ where $p$ runs over all prime numbers $\mathbb{P}$ and each $s_p \in \mathbb{N} \cup \{ \infty \}$. The set $\mathbb{S}$ of all supernatural numbers forms a multiplicative semigroup with multiplication defined by exponent addition and the multiplicative semigroup $\mon$ embeds in $\mathbb{S}$ via unique factorization. If $s \in \mathbb{S}$, the corresponding additive subgroup is
\[
P_s = \{ q \in \mathbb{Q}~|~s.q \in \mathbb{S} \} \]
In \cite{BZ} it is shown that isomorphism classes of additive subgroups of $\mathbb{Q}$ correspond to equivalence classes on $\mathbb{S}$ defined by
\[
s \sim s'~\text{iff}~\begin{cases} s_p = \infty \Leftrightarrow s'_p = \infty \\
s_p = s'_p~\text{for all but at most finitely many $p$} \end{cases} \]
That is, we have natural identifications between points of $\widehat{\wis{C}}$ and equivalence classes $[ \mathbb{S} ]$ of supernatural numbers. The connection between the two descriptions
\[
[ \mathbb{A}^f ] = \wis{pts}(\widehat{\wis{C}}) = [ \mathbb{S} ] \]
is induced via the identifications $\mathbb{S} = \widehat{\mathbb{Z}} / \widehat{\mathbb{Z}}^*$. In this note we will always use the description via supernatural numbers and will denote the isomorphism class of the point $H_s$ by the corresponding equivalence class $[s]$. Note also from \cite[2.2]{LBSieve} or \cite[\S 8]{CC2} that we can view $\wis{pts}(\widehat{\wis{C}})$ as the moduli space of noncommutative spaces corresponding to UHF-algebras.

\subsection{The subobject classifier $\mathbf{\Omega}$ :}  The {\em subobject-classifier} $\mathbf{\Omega}$ of the topos $\widehat{\wis{C}}$ coincides by with the set of all {\em sieves}, that is with all right ideals in the monoid $\mathbb{N}^{\times}_+$, see \cite[I.4]{MM}. Any sieve $\wis{S} \in \mathbf{\Omega}$ is of the form
\[
\wis{S} = \bigcup_{i \in I} n_i \mathbb{N}^{\times}_+ \]
for some (finite or infinite) set of positive integers $n_i$ which are incomparable with respect to the ordering by division. As a subobject-classifier in a topos, $\mathbf{\Omega}$ is a (complete) {\em Heyting algebra}.

Recall from \cite[I.8]{MM} that a {\em Heyting algebra} $\wis{H}$ is a partially ordered set under $\leq$, with unique minimal element $\wis{0}$ and unique maximal element $\wis{1}$ and equipped with $2$-ary operations $\wedge$ and $\vee$ and an additional $2$-ary operation $\rightarrow$ satisfying the following conditions for all $x,y,z \in \wis{H}$:
\begin{enumerate}
\item{$x \wedge (y \wedge z) = (x \wedge y) \wedge z = x \wedge y \wedge z$}
\item{$x \vee (y \vee z) = (x \vee y) \vee z = x \vee y \vee z$}
\item{$x \wedge y = y \wedge x$ and $x \vee y = y \vee x$}
\item{$x \wedge x = x$ and $x \vee x = x$}
\item{$\wis{1} \wedge x = x$ and $\wis{0} \vee x = x$}
\item{$x \wedge (y \vee x) = x = (x \wedge y) \vee x$}
\item{$x \wedge (y \vee z) = (x \wedge y) \vee (x \wedge z)$ and $x \vee ( y \wedge z) = (x \vee y) \wedge (x \vee z)$}
\item{$x \leq y$ iff $x = x \wedge y$ iff $y = x \vee y$}
\item{$(x \wedge y) \leq z$ iff $x \leq (y \rightarrow z)$}
\end{enumerate}
In an Heyting algebra one can also define the $1$-ary 'negation' operation $\neg$ by $\neg x = (x \rightarrow \wis{0})$ and one always has $x \leq \neg \neg x$. A Heyting algebra is called {\em Boolean} if for all $x \in \wis{H}$ we have $\neg \neg x = x$, or equivalently that $x \vee \neg x = \wis{1}$.

We will now make these operations explicit for the subobject-classifier $\mathbf{\Omega}$ of the presheaf topos $\widehat{\wis{C}}$ using \cite[Prop. I.8.5]{MM}. Assume $\wis{S}, \wis{T} \in \mathbf{\Omega}$ with
\[
\wis{S} = \bigcup_{i \in I} n_i \mathbb{N}^{\times}_+ \quad \text{and} \quad \wis{T} = \bigcup_{j \in J} m_j \mathbb{N}^{\times}_+ \]
then
\[
\wis{S} \vee \wis{T} = \wis{S} \cup \wis{T} = \bigcup_{i \in I}  n_i \mathbb{N}^{\times}_+  \cup \bigcup_{j \in J} m_j \mathbb{N}^{\times}_+ \]
\[
\wis{S} \wedge \wis{T} = \wis{S} \cap \wis{T} = \bigcup_{(i,j) \in I \times J} lcm(n_i,m_j) \mathbb{N}^{\times}_+ \]
\[
\wis{1} = \mathbb{N}^{\times}_+ \quad \text{and} \quad \wis{0} = \emptyset \]
\[
\wis{S} \rightarrow \wis{T} = \bigcup_{\forall i : lcm(e,n_i) \in \wis{T}} e \mathbb{N}^{\times}_+ \]
\[
\neg \wis{S} = \begin{cases} \wis{0} & \text{if $\wis{S} \not= \wis{0}$} \\
\wis{1} & \text{if $\wis{S} = \wis{0}$} \end{cases}  \]

\subsection{The Grothendieck topologies on $\wis{C}$} In addition to the Heyting algebra structure on $\mathbf{\Omega}$ there is also a right $\mathbb{N}^{\times}_+$-action $\odot$ on $\mathbf{\Omega}$.
\[
\wis{S} \odot n = n^{-1}.\wis{S} \cap \mathbb{N}^{\times}_+ = \bigcup_{i \in I} \frac{n_i}{gcd(n_i,n)} \mathbb{N}^{\times}_+ \]
Recall from \cite[III.2]{MM} that a {\em Grothendieck topology} on $\wis{C}$ is a subset $\mathcal{G} \subset \mathbf{\Omega}$ satisfying the following properties
\begin{enumerate}
\item{$\wis{1} \in \mathcal{G}$}
\item{(stability) if $\wis{S} \in \mathcal{G}$, then $\wis{S} \odot \mathbb{N}^{\times}_+ = \{ \wis{S} \odot n~|~n \in \mathbb{N}^{\times}_+ \} \subset \mathcal{G}$}
\item{(transitivity) if $\wis{S} \in \mathcal{G}$ and if $\wis{R} \in \mathbf{\Omega}$ such that $\wis{R} \odot \wis{S} = \{ \wis{R} \odot s~|~s \in \wis{S} \} \subset \mathcal{G}$, then $\wis{R} \in \mathcal{G}$.}
\end{enumerate}
Observe that it follows that if $\wis{S} \in \mathcal{G}$ and $\wis{S} \subset \wis{S'}$ in $\mathbf{\Omega}$, then also $\wis{S'} \in \mathcal{G}$.  The coarsest Grothendieck topology on $\wis{C}$, with $\mathcal{G}_{ch} = \{ \wis{1} \}$, is called the {\em chaotic topology} on the category $\wis{C}$. There are uncountable many different Grothendieck topologies on $\wis{C}$.

\begin{proposition} Let $P \subset \mathbb{P}$ be a subset of prime numbers, then
\[
\mathcal{G}_P = \{ \wis{S} \in \mathbf{\Omega}~|~\exists m : m \mathbb{N}^{\times}_+ \subset \wis{S}~\text{and all prime divisors of $m$ belong to $P$} \} \]
is the smallest Grothendieck topology on $\wis{C}$ containing the sieves $S(p) = p \mathbb{N}^{\times}_+$ for all $p \in P$.
\end{proposition}

\begin{proof} One uses the third property of Grothendieck topologies to deduce that $S(m) = m \mathbb{N}^{\times}_+$ must belong to $\mathcal{G}_P$ if all prime divisors of $m$ belong to $P$ using
\[
S(m) \odot S(p) = \begin{cases} S(m)~\text{ if $ p \nmid m$} \\ S(\frac{m}{p})~\text{ if $p \mid m$} \end{cases} \]
But then all sieves $\wis{S}$ of the required form belong to $\mathcal{G}_P$ and one verifies easily that these sieves indeed form a Grothendieck topology.
\end{proof}

But there are plenty of other Grothendieck topologies. If $\mathcal{G}$ is a Grothendieck topology on $\wis{C}$ and if $\cup_{i \in I} n_i \mathbb{N}^{\times}_+ \in \mathcal{G}$, then also all $\cup_{i \in I} p_i \mathbb{N}^{\times}_+ \in \mathcal{G}$ for all possible prime divisors $p_i$ of $n_i$. Hence, to $\mathcal{G}$ one can associate a set of subsets of prime numbers $\{ P_j : j \in J \}$ with the $P_j$ minimal with respect to the property that $\cup_{p \in P_j} p \mathbb{N}^{\times}_+ \in \mathcal{P}$. 

\section{Topos theoretic covers of the arithmetic site}

We have seen that the subobject classifier of $\widehat{\wis{C}}$ is $\mathbf{\Omega}$ with its natural right $\mathbb{N}^{\times}_+$-action. The terminal object of the topos $\widehat{\wis{C}}$ is given by the functor
\[
\wis{1}~:~\wis{C} \rTo \wis{Sets}  \qquad \bullet \mapsto \{ \ast \} \]
where $\{ \ast \}$ is the singleton with trivial right $\mathbb{N}^{\times}_+$-action. Open objects in $\widehat{\wis{C}}$ correspond to subobjects of $\wis{1}$, that is to right $\mathbb{N}^{\times}_+$-maps $\{ \ast \} \rTo \mathbf{\Omega}$ so its image must be a sieve fixed under the $\mathbb{N}^{\times}_+$-action. There are only two such sieves: $\emptyset = \wis{0}$ and $\mathbb{N}^{\times}_+ = \wis{1}$. That is $\wis{Opens}(\widehat{\wis{C}}) = \{ \wis{0},\wis{1} \}$, or equivalently, the topos $\widehat{\wis{C}}$ is two-valued.

\subsection{The Diaconescu cover} For every Grothendieck topos $\mathcal{E}$ there exists a locale $\wis{X}$ and an open surjective geometric morphism $\wis{Sh}(\wis{X}) \rTo \mathcal{E}$, \cite[IX.9 Thm 1]{MM}. The locale $\wis{X}$ is then called a {\em Diaconescu cover} of $\mathcal{E}$. In this section we will determine the Diaconescu covers of $\wis{Sh}(\wis{C},\mathcal{G})$ of the arithmetic site equipped with a Grothendieck topology.

As the topos of sheaves of sets on a site with underlying small category a partially ordered set is a locale by \cite[IX.5 Thm 1]{MM}, we first  introduce the poset corresponding to the monoid category $\wis{C}$, and call it the {\em big cell}.

\begin{definition} The {\em big cell} $\wis{D}$ is the category with one object $[k]$ for every strictly positive integer $k \in \mathbb{N}_+$ and with morphisms induced by the partial ordering by reverse division
\[
[k ] \leq [l] \quad \text{if and only if} \quad l | k \]
That is, there is exactly one arrow $[k] \rTo [l]$ if and only if $[k] \leq [l]$. There is a natural covariant functor $\wis{D} \rTo^{\pi} \wis{C}$ sending each object $[k]$ to $\bullet$ and every morphism $\xymatrix{[k] \ar[r] & [l]}$ to the endomorphism corresponding to $\tfrac{k}{l}$.
\end{definition}

The terminology 'big cell' s due to John Conway in \cite{Conway}. In order to understand groups commensurable with the modular group $\Gamma$ occurring in moonshine, he looked at (projective classes of) lattices commensurable with the standard lattice $< \vec{e}_1,\vec{e}_2 >$ and showed that they have a canonical form $L = < M \vec{e}_1 + \frac{g}{h} \vec{e}_2, \vec{e}_2 >$ with $M \in \mathbb{Q}_+$ and $\frac{g}{h} \in \mathbb{Q}/\mathbb{Z}$. For any other such lattice $L' = < N \vec{e}_1 + \frac{i}{j} \vec{e}_2, \vec{e}_2 >$ he introduced a {\em hyper-distance}
\[
\delta(L,L') = det(\alpha D_{LL'}) \in \mathbb{Z} \]
where $\alpha$ is the minimal strictly positive rational number such that $\alpha D_{LL'} \in M_2(\mathbb{Z})$ with
\[
D_{LL'} = \begin{bmatrix} M & \frac{g}{h} \\ 0 & 1 \end{bmatrix} . \begin{bmatrix} N & \frac{i}{j} \\  0 & 1 \end{bmatrix}^{-1} \]
The hyperdistance is symmetric and depends only on the projective classes of lattices. We can turn the set of all classes of lattices into a graph by connecting lattices having prime hyperdistance. Conway calls the resulting graph the {\em big picture} and proved that is the product, over all prime numbers $p$, of the subgraphs of lattices with hyperdistance $p$, which is itself an infinite $p+1$-valent tree, \cite{Conway}. The big cell is the subgraph consisting of the lattices $L_M$ with $M \in \mathbb{N}_+$, and is therefore the product of an infinite number of graphs of type $A_{\infty}$.

A Grothendieck topology $\mathcal{G}$ on $\wis{C}$ can be extended, in a canonical way, to a Grothendieck topology $\mathcal{G}_c$ on $\wis{D}$ by taking for each object $[k]$ in $\wis{D}$ the set of sieves $\mathcal{G}_c([k]) = \mathcal{G}$ where we mean for each sieve $\wis{S} \in \mathbf{\Omega}$ that
\[
\wis{S}([k]) = \{ \xymatrix{[l] \ar[r] & [k]}~:~\frac{l}{k} \in \wis{S} \}. \]
With respect to these Grothendieck topologies on $\wis{D}$ and $\wis{C}$ the functor $\pi$ satisfies the covering lifting property (clp)  as well as the property of preserving covers of \cite[p. 507]{MM}. It follows that the conditions of \cite[Prop. IX.8.1]{MM} are satisfied and we obtain the following result.

\begin{proposition} For any Grothendieck topology $\mathcal{G}$ on $\wis{C}$, the induced map from $\pi$ on the corresponding sheaf toposes
\[
f~:~\wis{Sh}(\wis{D},\mathcal{G}_c) \rTo \wis{Sh}(\wis{C},\mathcal{G}) \]
is an open, surjective geometric morphism. 

As $\wis{D}$ is a poset, the topos $\wis{Sh}(\wis{D},\mathcal{G}_c)$ is {\em localic}, that is, $\wis{Sh}(\wis{D},\mathcal{G}_c)$ is a Diaconescu cover of $\wis{Sh}(\wis{C},\mathcal{G})$.
\end{proposition}

Recall that a {\em point} of a topos is a geometric morphisms from $\wis{Sets}$ to the topos. As composition of geometric morphisms remain geometric, the above functor $f$ defines a map on the level of points
\[
\wis{pts}(\wis{D},\mathcal{G}_c) \rTo \wis{pts}(\wis{C},\mathcal{G}) \]
We will now investigate this map in case $\mathcal{G}$ is the chaotic topology, that is, $\wis{Sh}(\wis{C},\mathcal{G}) = \wis{PSh(\wis{C})} = \widehat{\wis{C}}$ in which case $\wis{Sh}(\wis{D},\mathcal{G}_c) = \wis{PSh}(\wis{D})$ which topos we will denote by $\widehat{\wis{D}}$.

\begin{proposition} The points of $\widehat{\wis{D}}$ correspond to subsets $S \subset \mathbb{N}^{\times}_+$ closed under taking divisors and least common multiples. As such they correspond naturally  to supernatural numbers: $\wis{pts}(\widehat{\wis{D}}) = \mathbb{S}$.
\end{proposition}

\begin{proof} Recall from \cite[VII,\S 5-6]{MM} that a point in the presheaf topos correspond to covariant functors
\[
P~:~\wis{D} \rTo \wis{Sets} \]
satisfying the following conditions
\begin{itemize}
\item{There is some $[k]$ such that $P([k]) \not= \emptyset$}
\item{Given elements $a \in P([k])$ and $b \in P([l])$, there is an $[m]$ with maps $\xymatrix{[k] & [m] \ar[l]_u \ar[r]^v & [l]}$ and an element $c \in P([m])$ such that $u(c)=a$ and $v(c)=b$.}
\end{itemize}
As $\wis{D}$ does not have parallel arrows, the third condition of \cite[p. 386]{MM} is vacuous. Note that for all $[k]$ we have that either $P([k])$ is the empty set or a singleton. For, suppose $a \not= b$ in $P([k])$, apply the second condition to $k=l$. As there is a unique map from $[m]$ to $[k]$ it must follow that $a=u(c)=b$.

By the first condition we have some $[k]$ such that $P([k])= \{ \ast_k \}$. If $l | k$ there is a unique map $\xymatrix{[k] \ar[r]^u & [l]}$ and therefore $u(\ast_k) \in P([l]) \not= \emptyset$. Hence the set
\[
S_P = \{ k \in \mathbb{N}^{\times}_+~|~P([k]) \not= \emptyset \} \]
is closed under division. It is also closed under taking least common multiples, for if $k,l \in S_P$ then by the second condition there is an $m$ such that $k | m$ and $l | m$ such that for the unique maps $u$ and $v$ from $[m]$ to resp. $[k]$ and $[l]$ we have
\[
u(\ast_m) = \ast_l \quad \text{and} \quad v(\ast_m) = \ast_l \]
as $lcm(k,l) | m$ it follows that also $lcm(k,l) \in S_P$.

To a subset $S \subset \mathbb{N}^{\times}_+$ satisfying these conditions one can associate a supernatural number $s \in \mathbb{S}$ with
\[
s = \prod_{p \in \mathbb{P}} p^{s_p} \]
where $s_p$ is the maximal exponent $d$ such that $p^d \in S$. Conversely, if $s \in \mathbb{S}$ then the subset
\[
S = \{ n \in \mathbb{N}^{\times}_+~|~n | s \} \]
is closed under division and least common multiples. 
\end{proof}

\subsection{The locale $\mathbb{S}$} The {\em terminal object} in $\widehat{\wis{D}}$ is given by the contravariant functor
\[
\wis{1}~:~\wis{D} \rTo \wis{Sets} \qquad \begin{cases} [k] \mapsto \{ \ast \} \\
\xymatrix{[k] \ar[r] & [l]} \mapsto id_{\{ \ast \}} \end{cases} \]
By definition, the set of opens in the topos $\widehat{\wis{D}}$ correspond to subfunctors of $\wis{1}$, which are described in the following result.

\begin{proposition} $\wis{Opens}(\widehat{\wis{D}})$ correspond to subsets of $\mathbb{N}^{\times}_+$ closed under taking multiples, that is, there is a natural one-to-one correspondence
\[
\wis{Opens}(\widehat{\wis{D}}) \leftrightarrow \mathbf{\Omega}  \]
between the opens in $\widehat{\wis{D}}$ and sieves in $\wis{C}$.
\end{proposition}

\begin{proof} A subfunctor $\wis{S} \in \wis{Sub}(\wis{1})$ assigns to each $[k]$ a set $S_k$ which is either $\{ \ast \}$ or the emptyset. If $k | l$ we must have the commuting diagram
\[
\xymatrix{\wis{1} : & \{ \ast \} \ar[r]^{id} & \{ \ast \} \\
\wis{S} : & S_k \ar[u] \ar[r]_{id | S_k} & S_l \ar[u]} \]
so the set of $k \in \mathbb{N}^{\times}_+$ such that $S_k = \{ \ast \}$ is closed under taking multiples, so is a right ideal of the monoid $\mathbb{N}^{\times}_+$, that is, a sieve in $\wis{C}$.
\end{proof}

With operations as before, $\mathbf{\Omega}$ is a {\em frame}, that is, a lattice with finite meets, finite and infinite joins, satisfying the infinite distributive law. More precisely, we have for all $\wis{S},\wis{T},\wis{T}_i \in \mathbf{\Omega}$:
\begin{enumerate}
\item{$\wis{S} \leq \wis{T}$ iff $\wis{S} \subset \wis{T}$.}
\item{$\wis{S} \wedge \wis{T} = \wis{S} \cap \wis{T}$.}
\item{$\wis{S} \vee \wis{T} = \wis{S} \cup \wis{T}$.}
\item{$\wis{S} \wedge (\vee_i \wis{T}_i) = \vee_i (\wis{S} \cap \wis{T}_i)$.}
\item{$\wis{0} = \emptyset$ and $\wis{1} = \mathbb{N}_+^{\times}$.}
\end{enumerate}

If $\wis{cat}$ is a poset category, then $\widehat{\wis{cat}}$ is generated by the subobjects of the terminal object, and is called a {\em localic topos} as $\widehat{\wis{cat}} \simeq \wis{Sh}(\wis{X})$, the category of sheaves on a {\em locale} $\wis{X}$, where the category of locales is defined to be the opposite of the category of frames, see for example \cite[Chp. IX]{MM}. The locale $\wis{X}$ has a frame isomorphic to the frame of open sets of a topological space $X$ if and only if $\wis{X}$ has {\em enough points}, that is, if elements of the corresponding frame can be distinguished by points of $\wis{X}$, \cite[Prop. IX.3.3]{MM}. Here, a point of $\wis{X}$ is a frame morphism
\[
p^{-1}~:~\mathcal{O}(X) = \wis{Open}(\widehat{\wis{cat}}) \rTo \{ \wis{0},\wis{1} \} \]

\begin{theorem} The localic topos $\widehat{\wis{D} }$ is equivalent to the category of sheaves of sets on the topological space $\mathbb{S}$
\[
\widehat{\wis{D} } = \wis{Sh}(\mathbb{S}) \]
of all supernatural numbers, with open sets corresponding to $\wis{S} = \cup_i n_i \mathbb{N}^{\times}_+ \in \mathbf{\Omega}$ where
\[
\mathbb{X}_l(\wis{S}) = \mathbb{X}_l(\cup_i n_i \mathbb{N}^{\times}_+) = \{ s \in \mathbb{S}~|~\exists i~:~n_i | s \} \]
\end{theorem}

\begin{proof} As $\wis{D}$ is a poset category, $\widehat{\wis{D} } \simeq \wis{Sh}(\wis{X})$ for some locale $\wis{X}$. A point of $\wis{X}$ is a frame morphism
\[
p^{-1}~:~\mathcal{O}(\wis{X}) = \wis{Open}(\widehat{\wis{D} }) = \mathbf{\Omega} \rTo \{ \wis{0},\wis{1} \} \]
or equivalent, a map satisfying the following conditions:
\begin{enumerate}
\item{$p^{-1}(\N^{\times}_+) = 1$}
\item{$p^{-1}(\wis{S} \wedge \wis{S'})=0$ iff $p^{-1}(\wis{S})=0$ or $p^{-1}(\wis{S'})=0$}
\item{$p^{-1}(\vee \wis{S}_i)=0$  iff $p^{-1}(\wis{S}_i ) = 0$  for all $i$}
\end{enumerate}
If $P_s$ is the point of $\widehat{\wis{D} }$ corresponding to the supernatural number $s \in \mathbb{S}$, that is, $P_s = \{ [k]~:~k | s \}$, then there is a  corresponding point $p_s$ of the locale $\wis{X}$ determined by the map
\[
p_s^{-1} = \mathbf{\Omega} \rTo \{ 0,1 \} \qquad p^{-1}(\cup_i n_i \N^{\times}_+) = \begin{cases} 0 & \text{iff $\forall i : n_i \nmid s$} \\ 1 & \text{iff $\exists i : n_i | s$} \end{cases} \]
as one easily verifies that this map satisfies the required properties.

But then, the locale $\wis{X}$ has enough points. For any two distinct elements $\wis{S},\wis{S}' \in \mathbf{\Omega} = \mathcal{O}(\wis{X})$, there is some point $p$ such that $p^{-1}(\wis{S}) \not= p^{-1}(\wis{S}')$. Indeed, let $k \in \wis{S} - \wis{S}'$ then $p_k^{-1}(\wis{S})=1$ whereas $p_k^{-1}(\wis{S}') = 0$.
To finish the proof we have to show that $\wis{pts}(\wis{X}) = \mathbb{S}$. 

Let $p^{-1} : \mathbf{\Omega} \rTo \{ \wis{0},\wis{1} \}$ be a point of $\wis{X}$ and consider the set $D_p = \{ n \in \mathbb{N}_+~|~p^{-1}(n \mathbb{N}^{\times}_+) = 1 \}$. As $n \mathbb{N}^{\times}_+ \wedge m \mathbb{N}^{\times}_+ = lcm(m,n) \mathbb{N}^{\times}_+$ it follows from the second condition on $p^{-1}$ that $D_p$ is closed under divisors and least common multiples, whence corresponds to a supernatural number $s \in \mathbb{S}$.
\end{proof}

This {\em localic topology} on $\mathbb{S}$ has a countable basis of open sets of the form
\[
\mathbb{X}_l(n) = \{ s \in \mathbb{S}~:~n | s \} \]
Observe that these opens are clopen in the compact Hausdorff topology on $\mathbb{S}$ via the identification with $\prod_{p \in \mathbb{P}} (\mathbb{N} \cup \{ \infty \})$ (see also the next section).

\vskip 3mm

\subsection{The projection map} We have already seen that the functor $\pi~:~\wis{D}  \rTo \wis{C}$ induces an open surjective geometric morphism
$\widehat{\wis{D} } \rOnto^f \widehat{\wis{C}}$. We now want to verify that the corresponding map between the points 
\[
\wis{pts}(\widehat{\wis{D} }) = \mathbb{S} \rTo \wis{pts}(\widehat{\wis{C}}) = [ \mathbb{S} ] \]
is the natural map, that is sending a supernatural number to its equivalence class. This requires a chase through the defining functors.

\vskip 3mm

We have identified points of $\widehat{\wis{D} }$ with the set of supernatural numbers $\mathbb{S}$. The point corresponding to $s \in \mathbb{S}$ is defined by the covariant flat functor
\[
P_s~:~\wis{D}  \rTo \wis{Sets} \qquad [k] \mapsto \begin{cases} \{ \ast_k \}~& \text{if $k | s$} \\ \emptyset~& \text{if $k \nmid s$} \end{cases} \]
\[
\xymatrix{[k] \ar[r] &  [l]} \mapsto \begin{cases} \{ \ast_k \} \rTo \{ \ast_l \}~& \text{if $l | k | s$} \\
\emptyset \rTo \{ \ast_k \}~& \text{if $l | s$ but $k \nmid s$} \\
\emptyset \rTo \emptyset~& \text{if $k,l \nmid s$} \end{cases} \]
This functor induces a geometric morphism $p(s)~:~\wis{Sets} \rTo \widehat{\wis{D} }$ determined by the adjoint functors
\[
\xymatrix{ \wis{Sets} \ar@/^2ex/[rr]^{p(s)_*} & & \widehat{\wis{D} } \ar@/^2ex/[ll]^{p(s)^*}} \quad p(s)_* = \underline{\wis{Hom}}_{\wis{D} }(P_s,-) \quad p(s)^*= - \otimes_{\wis{D} } P_s  \]
see \cite[VII.5 Thm 2, VII.2 Thm 1]{MM}. Here, $\underline{\wis{Hom}}_{\wis{D} }(P_s,S)$ is the contravariant functor (presheaf)
\[
\wis{D}  \rTo \wis{Sets} \qquad [k] \mapsto \begin{cases} \wis{Hom}_{\wis{Sets}}(\{ \ast_k \},S) = S~& \text{if $k | s$} \\ \wis{Hom}_{\wis{Sets}}(\emptyset, S) = \{ i_k \}~& \text{if $k \nmid s$} \end{cases} \]
and with $\xymatrix{[k] \ar[r] & [l]}$ mapping to the composition if $l | k | s$
\[
\xymatrix{ \{ \ast_k \} \ar@{.>}[r] \ar[d] & S \\ \{ \ast_l \} \ar[ru] &} \qquad
\xymatrix{ \emptyset \ar@{.>}[r]^{i_k} \ar[d] & S \\ \{ \ast_l \} \ar[ru] &} \qquad
\xymatrix{ \emptyset \ar@{.>}[r] \ar[d] & S \\ \emptyset \ar[ru]_{i_l} &}
 \]
and to the canonical maps $S \mapsto \{ i_k \}$ if $l | s$ but $k \nmid s$ and $\{ i_l \} \rTo \{ i_k \}$ if $k,l \nmid s$.

In general, a presheaf $G \in \widehat{\wis{D} }$ is a contravariant functor $G : \wis{D}  \rTo \wis{Sets}$ and is defined by assigning a set $G_k$ to all $[k]$, and a morphism $g_{k,l}~:~G_l \rTo G_k$ whenever $l | k$. Working through the definition of \cite[p. 356]{MM} we find that $G \otimes_{\wis{D} } P_s$ is the set of equivalence classes of the set
\[
\bigsqcup_{k | s} G_k \quad \text{with} \quad g_k \sim g_l \quad \text{iff} \quad g_{k,l}(g_l)=g_k \]

The functor $\pi~:~\wis{D} \rTo \wis{C}$ induces by \cite[VII.2 Thm 2]{MM} a geometric morphism
\[
\xymatrix{\widehat{\wis{D} } \ar@/^2ex/[rr]^{\pi_*} & & \hat{\wis{C}} \ar@/^2ex/[ll]^{\pi^*}} \]
where $\pi^*$ assigns to a $\mathbb{N}^{\times}_+$-set $M$, that is a functor $M : \wis{C}^{op} \rTo \wis{Sets}$, the composition
\[
\pi^*(M)~:~\wis{D} ^{op} \rTo^{\pi^{op}} \wis{C}^{op} \rTo^M \wis{Sets} \]
Concretely, $\pi^*(M)$ is the presheaf in $\widehat{\wis{D} }$ with all $G_k = M$ and the maps $g_{k,l}$ are given by the right action of $\frac{k}{l}$ on $M$ whenever $l | k$.

\begin{theorem} The functor $\pi~:~\wis{D}  \rTo \wis{C}$ induces a map
\[
\pi~:~\mathbb{S} = \wis{pts}(\widehat{\wis{D} }) \rTo \wis{pts}(\hat{\wis{C}}) = \mathbb{S}/ \sim \]
sending a supernatural number $s$ to its equivalence class $[s]$.
\end{theorem}

\begin{proof} In view of the foregoing, we have to prove that the composition $p(s)^* \circ \pi^*$
\[
\xymatrix{\wis{Sets} \ar@/^2ex/[rr]^{p(s)_*} & & \widehat{\wis{D} } \ar@/^2ex/[rr]^{\pi_*} \ar@/^2ex/[ll]^{p(s)^*}& & \hat{\wis{C}} \ar@/^2ex/[ll]^{\pi^*} \ar@/^8ex/[llll]^{- \otimes_{\wis{C}} \mathbb{Q}_+(s)}}
\]
is naturally isomorphic to the functor $- \otimes_{\wis{C}} \mathbb{Q}_+(s)$ where $\mathbb{Q}_+(s) = \{ \frac{n}{k}~|~n \in \mathbb{N}_+, k | s \}$ is the set with left $\mathbb{N}^{\times}_+$-action by multiplication corresponding to the flat functor $\wis{C} \rTo \wis{Sets}$ determining the point $[s] \in \wis{pts}(\hat{\wis{C}})$.
Let $M$ be a set with a right $\mathbb{N}^{\times}_+$-action, then $p(s)^*(\pi^*(M))$ is the set of equivalence classes of
\[
\bigsqcup_{k | s} M \quad \text{with} \quad [m']_k \sim [m]_1 \quad \text{iff} \quad m'=m.k \]
that is, we should view $m \in M$ in component $k$ as $m.\frac{1}{k}$ and hence we recover indeed the set $M \otimes_{\wis{C}} \mathbb{Q}_+(s)$.
\end{proof}

\subsection{The Barr cover} Barr's theorem \cite[Thm. IX.9.2]{MM} asserts that for every Grothendieck topos $\mathcal{E}$ there exists a complete Boolean algebra $\wis{B}$ and a surjective geometric morphism $\wis{Sh}(\wis{B}) \rOnto \mathcal{E}$. The purpose of this theorem is to generalize the Godement resolution by flabby sheaves to the context of cohomology of arbitrary toposes. We will now describe the Barr cover of the arithmetic site.

\begin{theorem} Consider the complete Boolean algebra
\[
\wis{B} = \wis{2}^{\mathbf{\Omega}} = \prod_{\wis{S} \in \mathbf{\Omega}} \{ \wis{0},\wis{1} \} \]
then there are surjective geometric morphisms
\[
\wis{Sh}(\wis{B}) \rOnto  \wis{Sh}(\mathbb{S}) \rOnto \widehat{\wis{C}} \]
That is, $\wis{B}$ is the Barr cover of the arithmetic site.
\end{theorem}

\begin{proof} We have already seen that the locale $\mathbb{S}$ is a topological space with the frame of opens $\wis{Open}(\mathbb{S}) \leftrightarrow \mathbf{\Omega}$. This allows us to determine the Heyting algebra structure on $\wis{Open}(\mathbb{S})$ as in \cite[I.8]{MM} and verify that it coincides with the Heyting algebra structure on $\mathbf{\Omega}$ given before. In particular we have for each non-empty $\wis{S}$ that its negation $\neg \wis{S} = \wis{0} = \emptyset$.

Following the construction of the proof of  \cite[Lemma IX.9.3]{MM}, we have to determine the frame of $\neg \neg$-fixed points of the frame $\wis{Open}(\mathbb{S} - \wis{S})$ of the closed sublocale $\mathbb{S} - \wis{S}$ of $\mathbb{S}$. Now, by definition
\[
\wis{Open}(\mathbb{S} - \wis{S}) = \{ \wis{T} \in \mathbf{\Omega}~|~\wis{S} \leq \wis{T} \} \cup \{ \wis{0} \} \]
and by the above remark on the negation in $\wis{Open}(\wis{S})$ we have that the double-negation frame
\[
\wis{Open}(\mathbb{S} - \wis{S})_{\neg \neg} \simeq \{ \wis{0}, \wis{1} \} \]
The claims now follow from the proof of \cite[Lemma IX.9.3]{MM}.
\end{proof}

\section{Noncommutative covers of the arithmetic site}

Rather than studying Grothendieck topologies on $\wis{C}$ and their associated toposes of sheaves, we will consider in this section a number of ordinary topologies on the set of points $\wis{pts}(\widehat{\wis{C}}) = [ \mathbb{A}^f ] = [ \mathbb{S} ]$. We then construct a noncommutative topological space, in the sense of \cite{Bauer}, with commutative shadow $[ \mathbb{S} ]$ and a corresponding noncommutative Heyting algebra $\mathbf{\Theta}$, in the sense of \cite{Cvetko}, with commutative quotient the subobject classifier $\mathbf{\Omega}$.

These results are inspired by a remark of A. Connes in \cite[\S 5]{ConnesRH} where he asserts that the arithmetic site is but a semiclassical shadow of a still mysterious structure dealing with compactifications of $\wis{Spec}(\mathbb{Z})$. For this reason we first consider three topologies on the points $[ \mathbb{S} ]$ of the arithmetic site and the corresponding continuous maps from $\wis{Spec}(\mathbb{Z})$. We then construct a noncommutative frame $\mathbf{\Theta}$, elements of which can be viewed as constructible truth fluctuations on open sets in the sieve topology on the finite ad\`ele classes.

\subsection{The arithmetic topology: } The locally compact topology on the finite ad\`eles $\mathbb{A}^f$ induces a weak remnant on the set of points of the arithmetic site, which we call the {\em arithmetic topology}. If we identify $\mathbb{S} = \prod_{p \in \mathbb{P}} (\mathbb{N} \cup \{ \infty \}$, view each of the factors as a point point compactification of the discrete topology on $\mathbb{N}$ then the arithmetic topology is the product topology turning $\mathbb{S}$ into a compact Hausdorff. The arithmetic topology is the induced topology on the equivalence classes $[ \mathbb{S} ]$ and was described in \cite[Thm. 1]{LBSieve}. It has a countable basis of open sets, corresponding to finite subsets $P = \{ p_1,\hdots,p_k \}$ of prime numbers, defined by
\[
\mathbb{X}_a(P) = \{ [s] \in [\mathbb{S} ]~|~\forall p \in P~:~s_p \not= \infty \} \]
(note that this does not depend on the choice of representative in $[s]$).
There is an obvious connection with the Zariski topology on the prime spectrum $\wis{Spec}(\mathbb{Z})$ which is $\mathbb{P} \cup \{ 0 \}$ with open sets of the form $\mathbb{X}_Z(n) = \{ p \in \mathbb{P}~|~p \nmid n \} \cup \{ 0 \}$. We can map a prime number $p$ to the class determined by the element  $(1,\hdots,1,0,1,\hdots) \in \widehat{\mathbb{Z}}$ with a zero in the factor $\widehat{\mathbb{Z}}_p$ and ones in the other factors. A direct calculation proves the following result.

\begin{proposition} The inclusion map 
\[
\wis{Spec}(\mathbb{Z}) \rInto^i \wis{pts}(\widehat{\wis{C}}) = [ \mathbb{S} ] \quad \text{given by} \quad p \mapsto [ p^{\infty} ] \]
is continuous with respect to the Zariski topology on $\wis{Spec}(\mathbb{Z})$ and the arithmetic topology on $[ \mathbb{S} ]$ as $i^{-1}(\mathbb{X}_a(P)) = \mathbb{X}_Z(\prod_{p \in P} p)$. As a consequence, the direct image sheaf $i^* \mathcal{O}_{\mathbb{Z}}$ has as sections over the arithmetic topology on $[\mathbb{S} ]$
\[
\Gamma(i^* \mathcal{O}_{\mathbb{Z}}, \mathbb{X}_a(P)) = \mathbb{Z}[\frac{1}{\prod_{p \in P} p}] \]
The stalk of $i^* \mathcal{O}_{\mathbb{Z}}$ in a point $[s]$ of the arithmetic site is then
\[
(i^* \mathcal{O}_{\mathbb{Z}})_{[s]} = \{ \frac{a}{b}~|~(b,p)=1~\forall p \in \mathbb{P} : s_p = \infty \} \]
so, in particular, the stalk in $[p^{\infty}]$ is the local ring $\mathbb{Z}_p$.
\end{proposition}

\subsection{The sieve topology: } In \cite{LBSieve} we introduced another topology on the set of points of the arithmetic site, called the {\em sieve topology} as the basic open sets correspond to the sieves $\wis{S} \in \mathbf{\Omega}$ and are defined as
\[
\mathbb{X}_s(\wis{S}) = \{ [ s ] \in [ \mathbb{S} ]~|~\exists n_1,n_2,\hdots \in \wis{S}~:~\prod_{i=1}^{\infty} n_i | s \} \]
(this definition does not depend on the choice of representative in $[s]$). The 'morale' behind this definition is that $\mathbb{X}_s(\wis{S})$ are precisely the points of the arithmetic site which are also points in the monoid category of $\wis{S} \cup \{ id \}$. In \cite[Thm. 6]{LBSieve} it was shown that this sieve topology has some properties one might expect of the mythical space $\overline{\wis{Spec}(\mathbb{Z})}$: it is compact, does not admit a countable basis of opens, every non-empty open set is dense and it satisfies the $T_1$-separation property on incomparable points. Further note that the frame of open sets in the sieve topology is isomorphic to the lattice structure of the Heyting algebra $\mathbf{\Omega}$.

For the connection with the arithmetic topology, note that if $P = \{ p_1,\hdots,p_k \}$ we have
\[
\mathbb{X}_a(P) = [ \mathbb{S} ] - \mathbb{X}_s( p_1 \mathbb{N}^{\times}_+ \cup \hdots \cup p_k \mathbb{N}^{\times}_+) \]
That is, basic opens of the arithmetic topology are closed in the sieve topology. Still, we have a nice connection with $\wis{Spec}(\mathbb{Z})$. This time we map the prime ideal $(p)$ to the class of the element $(0,\hdots,0,1,0,\hdots) \in \widehat{\mathbb{Z}}$ with a one in the factor corresponding to $\widehat{\mathbb{Z}}_p$ and zeroes elsewhere. In this case we calculate:

\begin{proposition} The inclusion map
\[
\wis{Spec}(\mathbb{Z}) \rInto^j \wis{pts}(\widehat{\wis{C}}) = [ \mathbb{S} ] \quad \text{given by} \quad p \mapsto [ \prod_{q \not= p} q^{\infty} ] \]
is continuous with respect to the Zariski topology on $\wis{Spec}(\mathbb{Z})$ and the sieve topology on $[ \mathbb{S} ]$ as for each sieve $\wis{S} = \cup_{i \in I} n_i \mathbb{N}^{\times}_+$ we have with $gcd(\wis{S}) = gcd(n_i, i \in I)$
\[
j^{-1}(\mathbb{X}_s(\wis{S})) = \mathbb{X}_a(gcd(\wis{S})) \]
As a consequence the direct image sheaf $j^* \mathcal{O}_{\mathbb{Z}}$ has as sections with respect to the sieve topology
\[
\Gamma(j^* \mathcal{O}_{\mathbb{Z}},\mathbb{X}_s(\wis{S})) = \mathbb{Z}[\frac{1}{gcd(\wis{S})}] \]
Therefore, the stalk of $j^* \mathcal{O}_{\mathbb{Z}}$ in a point $[s]$ is equal to
\[
(j^* \mathcal{O}_{\mathbb{Z}})_{[s]} = \mathbb{Z}[ \frac{1}{p}~:~p^{\infty} | s ] \]
or, in terms of the additive subgroups of $\mathbb{Q}$
\[
(j^* \mathcal{O}_{\mathbb{Z}})_{[s]} = \bigcap_{[t]=[s]} P_t \]
In particular, the stalk of $j^* \mathcal{O}_{\mathbb{Z}}$ at the image $j(p)$ is the localization $\mathbb{Z}_p$.
\end{proposition}

This result shows that the sieve topology and its induced structure sheaf from the connection with $\wis{Spec}(\mathbb{Z})$ is close to the intended structure sheaf of the arithmetic site as defined by Connes and Consani in \cite{CC} and \cite{CC2}. The {\em stalk} of the structure sheaf $\mathcal{O}$ at a point $[s] \in \wis{pts}(\widehat{\wis{C}})$, which as we have seen corresponds to the (isomorphism class) of the additive subgroup $P_s$ of $\mathbb{Q}$, will be the semiring $P_{s,max} = (P_s \cup \{ - \infty \},max,+)$ as sub-semiring of $\mathbb{Q}_{max}$.

Clearly, if one ever wants to have a concrete realization of this "stalk" as a genuine stalk over some topology on the set of all points of the arithmetic site one has to dispose of the ambiguity that different semirings $P_{s,max}$ defining the same point $[s]$ are only isomorphic as semirings, and not equal. Replacing $P_{s.max}$ by the intersection $\cap_{[t]=[s]} P_{t,max}$ is an obvious choice, which gives us the exact stalk of the sheaf $j^* \mathcal{O}_{\mathbb{Z}}$ in the sieve topology.

\subsection{The patch topology: } At this moment we have two topologies on $\wis{pts}(\widehat{\wis{C}}) = [ \mathbb{S} ]$ and two natural continuous maps from $\wis{Spec}(\mathbb{Z})$ to $[ \mathbb{S} ]$ such that in both cases at least the stalks of the direct image structure sheaf give the expected ring $\mathbb{Z}_p$. The arithmetic topology might be better to transport arithmetic (ad\`elic) information to the arithmetic site whereas the sieve topology, by its very nature, is better at encoding topos-theoretic information.

In order to have the best of both worlds, we consider a common refinement of the two topologies, the {\em patch topology} (or {\em constructible topology}) with respect to the sieve topology. That is, open sets in the patch topology are unions of locally closed subsets in the sieve topology on $\wis{pts}(\widehat{\wis{C}})$. That is, a basic open set in the patch topology is of the form
\[
\mathbb{X}_p(\wis{S},\wis{T}) = \mathbb{X}_s(\wis{S}) \cap \mathbb{V}_s(\wis{T}) \]
where $\wis{S},\wis{T} \in \mathbf{\Omega}$ and $\mathbb{V}_s(\wis{T}) = [ \mathbb{S} ] - \mathbb{X}_s(\wis{T})$. As we have seen that the basic open sets in the arithmetic topology are closed in the sieve topology we have:

\begin{proposition} The patch topology on the points of the arithmetic site is a common refinement of both the arithmetical and the sieve topology. It is the topology generated by both the open and closed sets of the sieve topology.
\end{proposition}

\subsection{A noncommutative topology with shadow $\mathbf{\Omega}$}

Our aim is to equip the set of points $\mathbb{E} = [ \mathbb{S} ] \times \{ 0,1 \}$ with a {\em noncommutative topology} such that its classical commutative shadow gives us  the sieve topology on $[ \mathbb{S} ]$. Here, a noncommutative topology has a corresponding 'frame' of open sets $\mathbf{\Theta}$ which will be a skew-lattice on which the operations $\wedge$ and $\vee$ are no longer commutative. By the shadow property we will mean that any skew-lattice map to a (commutative) lattice must factor through the lattice $\mathbf{\Omega}$ of all opens of the sieve topology on $[\mathbb{S}]$. 

It turns out, see for example \cite{Bauer} or \cite{Cvetko}, that in passing to the noncommutative setting we need to sacrifice either the top element $\wis{1}$ or the bottom element $\wis{0}$ in the definition of a skew-lattice. In this section we choose to sacrifice $\wis{1}$ as our construction will involve taking local sections of sheaves and whereas there is always a unique such section over $\wis{0} = \emptyset$ there are usually lots of global sections.

\begin{definition} A {\em noncommutative frame} is a set $\mathbf{\Theta}$ with $2$-ary operations $\wedge$ and $\vee$ (intersection resp. union of open sets) and a minimal element $\wis{0}$ such that for all $x,y,z \in \mathbf{\Theta}$ the following conditions are satisfied:
\begin{enumerate}
\item{$x \wedge (y \wedge z) = (x \wedge y) \wedge z = x \wedge y \wedge z$}
\item{$x \vee ( y \vee z) = (x \vee y) \vee z = x \vee y \vee z$}
\item{$x \wedge x = x$}
\item{$x \vee x = x$}
\item{$x \wedge ( x \vee y) = x = x \vee (x \wedge y)$}
\item{$(y \vee x) \wedge x = x = (y \wedge x) \vee x$}
\item{$x \wedge \wis{0} = \wis{0} = \wis{0} \wedge x$}
\item{$x \vee \wis{0} = x = \wis{0} \vee x$}
\item{$x \wedge (y \vee z) = (x \wedge y) \vee (x \wedge z)$}
\item{$(y \vee z) \wedge x = (y \wedge x) \vee (z \wedge x)$}
\item{$x \wedge y \wedge x = x \wedge y$}
\item{$x \vee y \vee x = y \vee x$}
\item{$x \wedge y \wedge z = x \wedge z \wedge y$}
\end{enumerate}
Here, (1)-(4) state that the operations are associative and idempotent, (5)-(6) are called the absorption properties, (9)-(10) are the strongly distributive properties and (11)-(12) are the left-handed properties, see \cite{Bauer} for more details. $\mathbf{\Theta}$ becomes a partially ordered set via
\[
x \leq y \Leftrightarrow x \wedge y = x = y \wedge x \Leftrightarrow x \vee y = y = y \vee x \]
and (7)-(8) state that $\wis{0}$ is a minimal element with respect to this partial order. Condition (13) follows from the others and states that $(\mathbf{\Theta},\wedge)$ is a left-normal semigroup.
\end{definition}

One defines an equivalence relation on the elements of a noncommutative frame $\mathbf{\Theta}$ by
\[
x \sim y \Leftrightarrow x \wedge y \wedge x = x~\text{and}~y \wedge x \wedge x = y \]
A result of Jonathan Leech, see for example \cite[Thm. 2.1]{Bauer} asserts that the induced operations turn the equivalence classes $[ \mathbf{\Theta} ]$ into a (commutative) lattice without top element, and moreover, any noncommutative frame morphism $\mathbf{\Theta} \rTo \mathbf{\Lambda}$ where $\mathbf{\Lambda}$ is a lattice factors through $[ \mathbf{\Theta} ]$.

The main result of Bauer et al. \cite{Bauer} on noncommutative Priestly duality tells us that any noncommutative topology $\mathbf{\Theta}$ with commutative semiclassical shadow $[ \mathbf{\Theta} ] \simeq \mathbf{\Omega}$ comes from the local sections of a sheaf with respect to the patch topology for $\mathbf{\Omega}$. 

We will now construct the simplest such sheaf in the case of interest to us, that is, when $\mathbf{\Omega}$ is the commutative frame of opens of the sieve topology on the points of the arithmetic site. 

\begin{definition} The sheaf $\mathcal{O}^c$ of {\em constructible truth fluctuations} has as sections over the open set $\mathbb{X}_s(\wis{S})$ for $\wis{S} \in \mathbf{\Omega}$
\[
\Gamma(\mathcal{O}^c,\mathbb{X}_s(\wis{S})) = \{ x : \mathbb{X}_s(\wis{S}) \rTo \mathbb{B} = \{ 0,1 \}~\text{continuous} \} \]
Here, continuous is with respect to the discrete topology on the Boolean semifield $\mathbb{B}= \{ 0,1 \}$ and the induced patch topology on $\mathbb{X}_s(\wis{S})$. Alternatively, one can think of $x$ as the characteristic function of an open subset of $\mathbb{X}_s(\wis{X})$ in the patch topology.
\end{definition}

Using the sheaf of constructible truth fluctuations $\mathcal{O}^c$ we will now construct the noncommutative frame $\mathbf{\Theta}$ which is the set of 'opens' of a noncommutative topology on $\mathbb{E} = [ \mathbb{S} ] \times \mathbb{B}$.

\begin{theorem} \label{noncframe} Let $\mathbf{\Theta}$ be the set of all pairs $(\wis{S},x)$ where $\wis{S} \in \mathbf{\Omega}$ and $x \in \Gamma(\mathbb{X}_s(\wis{S}),\mathcal{O}^c)$. On $\mathbf{\Theta}$ we define $2$-ary operations $\wedge$ and $\vee$ as follows:
\[
(\wis{S},x) \wedge (\wis{T},y) = (\wis{S} \wedge \wis{T}, x |_{\mathbb{X}_s(\wis{S}) \cap \mathbb{X}_s(\wis{T})}) \]
\[
(\wis{S},x) \vee (\wis{T},y) = (\wis{S} \vee \wis{T}, y |_{\mathbb{X}_s(\wis{T})} \cup x |_{\mathbb{X}_s(\wis{S}) - \mathbb{X}_s(\wis{T})}) \]
With respect to these operations, $\mathbf{\Theta}$ is a noncommutative frame with semiclassical shadow $[ \mathbf{\Theta} ] \simeq \mathbf{\Omega}$. there is an embedding of noncommutative frames $\mathbf{\Omega} \rInto \mathbf{\Theta}$ given by sending $\wis{S}$ to $(\wis{S},1 |_{\mathbb{X}_s(\wis{S})})$. Further, $(\wis{S},x)$ commutes with $(\wis{T},y)$ for $\wedge$ or $\vee$ if and only if $x |_{\mathbb{X}_s(\wis{S}) \cap \mathbb{X}_s(\wis{T})} = y |_{\mathbb{X}_s(\wis{S}) \cap \mathbb{X}_s(\wis{T})}$.

To each $(\wis{S},x) \in \mathbf{\Theta}$ we associate the 'open' set in $\mathbb{E} = [ \mathbb{S} ] \times \mathbb{B}$
\[
\mathbb{X}_{nc}(\wis{S},x) = \{ ([s],x([s]))~:~[s] \in \mathbb{X}_s(\wis{S}) \} \]
The natural projection om the first factor makes $\mathbb{E}$ into an \'etale space over $[ \mathbb{S} ]$ with respect to the patch topology on $[ \mathbb{S} ]$.
\end{theorem}

\begin{proof} All claims follow from computations, or by invoking \cite{Bauer}. Here, we will merely make the partial ordering explicit as well as the equivalence relation. Let $(\wis{S},x)$ and $(\wis{T},y)$ be elements of $\mathbf{\Theta}$. Then, if
\[
(\wis{S},x) \wedge (\wis{T},y) = (\wis{S} \wedge \wis{T}, x |_{\mathbb{X}_s(\wis{S}) \cap \mathbb{X}_s(\wis{T})}) = (\wis{S},x) \]
this implies that $\wis{S} \leq \wis{T}$. The reverse condition
\[
(\wis{T},y) \wedge (\wis{S},x) = (\wis{T} \wedge \wis{S}, y |_{\mathbb{X}_s(\wis{T}) \cap \mathbb{X}_s(\wis{S})}) = (\wis{S},x) \]
asserts that $\wis{S} \leq \wis{T}$ and that $y |_{\mathbb{X}_s(\wis{S})} = x$. That is, the partial ordering on $\mathbf{\Theta}$
\[
(\wis{S},x) \leq (\wis{T},y)~\Leftrightarrow~\wis{S} \leq \wis{T}~\text{and}~y |_{\mathbb{X}_s(\wis{S})} = x \]
is given by extension of local sections. Using this, we also have that if
\[
(\wis{S},x) \wedge (\wis{T},y) \wedge (\wis{S},x) = (\wis{S} \wedge \wis{T}, x |_{\mathbb{X}_s(\wis{S}) \cap \mathbb{X}_s(\wis{T})}) = (\wis{S},x) \]
is equivalent to $\wis{S} \leq \wis{T}$ whereas the other condition
\[
(\wis{T},y) \wedge (\wis{S},x) \wedge (\wis{T},y) = (\wis{S} \wedge \wis{T}, y |_{\mathbb{X}_s(\wis{T}) \cap \mathbb{X}_s(\wis{S})}) = (\wis{T},y) \]
is equivalent to $\wis{T} \leq \wis{S}$. That is, $(\wis{S},x) \sim (\wis{T},y)$ iff $\wis{S}=\wis{T}$. As a consequence, $[ \mathbf{\Theta} ]$ is just projection on the first factor so is isomorphic to the lattice $\mathbf{\Omega}$.
\end{proof}

We can, of course, repeat this construction starting from more involved sheaves with respect to the sieve topology on $[ \mathbb{S} ]$. A natural choice would be to consider $\mathbb{F}$ to be the disjoint union of all the stalks of the structure sheaf $j^* \mathcal{O}_{\mathbb{Z}}$. That is, $\mathbb{F}$ is equal to the following subset of $[\mathbb{S}] \times \mathbb{Q}_{max}$:
\[
\bigsqcup_{[s] \in \mathbb{S}} (j^* \mathcal{O}_{\mathbb{Z}})_{[s]} = \{ ([s],q) \in [\mathbb{S}] \times \mathbb{Q}_{max}~|~\forall [t]=[s]~:~q.t \in \mathbb{S} \} \cup \{ ([s],-\infty)~|~[s] \in \mathbb{S} \} \]
One can then consider the sheaf of {\em constructible functions} $\mathcal{O}$ having as its sections over $\mathbb{X}_s(\wis{S})$
\[
 \{ x : \mathbb{X}_s(\wis{S}) \rTo \mathbb{Q}_{max}~:~Graf(x) \subset \mathbb{F}~\text{and locally constant wrt. patch topology} \} \]
 Using the same operations as above we can then define a noncommutative frame structure on the collection of all $(\wis{S},x)$ with $\wis{S} \in \mathbf{\Omega}$ and $x \in \Gamma(\mathcal{O},\mathbb{X}_s(\wis{S}))$ having as its semiclassical shadow $\mathbf{\Omega}$.

\subsection{A noncommutative Heyting algebra with shadow $\mathbf{\Omega}$}

In \cite{Cvetko} K. Cvetko-Vah introduced the notion of a {\em skew-Heyting algebra}. The starting point is again a skew-lattice but, in view of the crucial role of the top $\wis{1}$ as truth in logic, this time with a top element but without bottom element, to ensure noncommutativity of $\wedge$ and $\vee$. The additional requirements for a skew-lattice with $\wis{1}$ to a a skew-Heyting algebra are by \cite[Thm. 3.2]{Cvetko} that for all $x,y,z$ and $u$ we have
\begin{enumerate}
\item{$(x \rightarrow y ) = ( y \vee x \vee y \rightarrow y)$}
\item{$(x \rightarrow x) = \wis{1}$}
\item{$(x \wedge (x \rightarrow y) \wedge x) = x \wedge y \wedge x$}
\item{$(y \wedge (x \rightarrow y)) = y$ and $((x \rightarrow y) \wedge y) = y$}
\item{$(x \rightarrow u \vee (y \wedge z) \vee u) = (x \rightarrow (u \vee y \vee u)) \wedge (x \rightarrow (u \vee z \vee u))$}
\end{enumerate}
If we exchange the definition of $\vee$ with that of $\wedge$ in theorem~\ref{noncframe}, reverse the partial ordering and define $\wis{1}=\emptyset$ one can extend the resulting skew-lattice with top to a skew-Heyting algebra as this is an instance of the partial maps with poset domain example of \cite[4.4]{Cvetko}. However, this is rather artificial and we want to keep the definitions of $\vee$ and $\wedge$ of theorem~\ref{noncframe} on $\mathbf{\Theta}$. For this reason we proceed differently and immediately define the $2$-ary operation $\rightarrow$ on $\mathbf{\Theta}$ (using the corresponding operation on $\mathbf{\Omega}$) as follows
\[
(\wis{S},x) \rightarrow (\wis{T},y) =  (\wis{S} \rightarrow \wis{T}, y \cup 1_{\mathbb{X}_s(\wis{S} \rightarrow \wis{T}) - \mathbb{X}_s(\wis{T})})  \]
Observe that in the Heyting algebra $\mathbf{\Omega}$ we always have that $\wis{T} \leq (\wis{S} \rightarrow \wis{T})$ so the extension of the function $y$ to $\mathbb{X}_s(\wis{S} \rightarrow \wis{T})$ makes sense.
Of course, as we do not have a top element in $\mathbf{\Theta}$ we cannot hope to satisfy condition $(2)$ of a skew-Heyting algebra. That is, in the noncommutative world we have to sacrifice absolute truth.

Based on the example below, K. Cvetko-Vak \cite{KarinNHA} formalized a {\em non-commutative Heyting algebra} to be  a noncommutative frame (with bottom $\wis{0}$)  with a $2$-ary operation $\rightarrow$ satisfying the following axioms. Here, $t$ is a distinguished element of $T$, the set of all top class elements of the poset.

(NH1) $x \rightarrow y = ( y \vee ( t \wedge x \wedge t) \vee y) \rightarrow y$,

(NH2) $x \rightarrow x = x \vee t \vee x$,

(NH3) $x \wedge (x \rightarrow y) \wedge x = x \wedge y \wedge x$,

(NH4) $y \wedge (x \rightarrow y) = y$ and $(x \rightarrow y) \wedge y = y$,

(NH5) $x \rightarrow (t \wedge (y \wedge z) \wedge t) = (x \rightarrow (t \wedge y \wedge t)) \wedge (x \rightarrow (t \wedge z \wedge t))$.

\begin{theorem} With the above definition of $(\wis{S},x) \rightarrow (\wis{T},y)$, the noncommutative frame $\mathbf{\Theta}$ becomes a noncommutative Heyting algebra. Moreover, any morphism of noncommutative Heyting algebras $\mathbf{\Theta} \rTo \mathbf{\Lambda}$ where $\Lambda$ is a Heyting algebra, factors through the Heyting algebra $\mathbf{\Omega}$. 
\end{theorem}

\begin{proof} Let $x=(\wis{S},x), y=(\wis{T},y), z=(\wis{V},z)$ and $t = (\wis{1},t)$. Using the definitions
\[
x \wedge y = (\wis{S} \wedge \wis{T},x),~\quad~ x \vee y = (\wis{S} \vee \wis{T}, y \cup x |_{\mathbb{X}_s(\wis{S})-\mathbb{X}_s(\wis{T})},~\quad~x \rightarrow y = (\wis{S} \rightarrow \wis{T},y \cup 1 |_{\mathbb{X}_s(\wis{S} \rightarrow \wis{T})- \mathbb{X}_s(\wis{T})} \]
one verifies that both sides of the equalities in (NH1)-(NH5) are equal, see also \cite[Example 3.1]{KarinNHA}.
\end{proof}

{\bf Acknowledgement: } I like to thank K. Cvetko-Vah for pointing out that the example above did not satisfy the axioms given in a previous version of this paper and for providing the corrected axioms (NH1)-(NH5).

\end{document}